\documentclass[12pt,draft]{amsart}

\usepackage[cp1251]{inputenc}
\usepackage[T2A]{fontenc}
\usepackage{amsmath}
\usepackage{amsfonts}
\usepackage{amssymb}
\usepackage{amsthm}
\usepackage{euscript}
\usepackage[dvips]{graphicx}
\usepackage[all]{xy}
\usepackage{delarray}
\usepackage{amscd}
\usepackage[ukrainian,russian,english]{babel}
\theoremstyle{plain}
\newtheorem{theoremEng}{\hspace{17pt}Theorem}
\newtheorem{lemmaEng}[theoremEng]{\hspace{17pt}Lemma}
\newtheorem{propositionEng}[theoremEng]{\hspace{17pt}Proposition}

\newtheorem{corollaryEng}[theoremEng]{\hspace{17pt}Corollary}
\theoremstyle{definition}

\newtheorem{remarkEng}[theoremEng]{\hspace{17pt}Remark}
\theoremstyle{plain}
\newtheorem*{theoremEng*}{\hspace{17pt}Theorem}
\newtheorem*{lemmaEng*}{\hspace{17pt}Lemma}
\newtheorem*{propositionEng*}{\hspace{17pt}Proposition}
\newtheorem*{statementEng*}{\hspace{17pt}Statement}
\newtheorem*{corollaryEng*}{\hspace{17pt}Corollary}
\theoremstyle{definition}
\newtheorem*{definitionEng*}{\hspace{17pt}Definition}
\theoremstyle{remark}
\newtheorem*{notationEng*}{\hspace{17pt}Notation}
\newtheorem*{remarkEng*}{\hspace{17pt}Remark}
\large \tolerance=4000

\title[Singular perturbed periodic operators]
{Singularly perturbed periodic and semiperiodic differential
operators}
\address{Institute of Mathematics NAS of Ukraine \\
    Tereshchenkivska str., 3 \\
    Kyiv\\
    Ukraine\\
    01601}
\author[V.A. Mikhailets, V.M. Molyboga]
       {V.A. Mikhailets, V.M. Molyboga}
\email[Vladimir A. Mikhailets]{mikhailets@imath.kiev.ua}
\email[Volodymyr M. Molyboga]{molyboga@imath.kiev.ua} \keywords{}
\subjclass[2000]{}

\begin{document}
\begin{abstract}
Qualitative and spectral properties of the form-sums
\begin{equation*}
  S_{\pm}(V):=D_{\pm}^{2m}\dotplus V(x),\quad m\in \mathbb{N},
\end{equation*}
in the Hilbert space $L_{2}(0,1)$ are studied. Here the periodic
$(D_{+})$ and the semiperiodic $(D_{-})$ differential operators
are $D_{\pm}:\,u\mapsto -i u'$, and $V(x)$ is a 1-periodic
complex-valued distribution in the Sobolev spaces
$H_{per}^{-m\alpha}$, $\alpha\in [0,1]$.
\end{abstract}
\maketitle
\section{Introduction and statement of results}\label{Int}
In this paper, we study the operators $S_{+}(V)$ and $S_{-}(V)$
that are not selfadjoint in general and given on the Hilbert space
$L_{2}(0,1)$ by two-terms differential expressions of an even
order, with a 1-periodic complex-valued potential $V(x)$, which is
a distribution in $\mathcal{D}'_{1}$, and periodic and
semiperiodic boundary conditions,
\begin{align*}
  L_{2}(0,1):\hspace{40pt} & S_{\pm}u\equiv S_{\pm}(V)u:=D_{\pm}^{2m}u+V(x)u, \\
  m\in \mathbb{N}\hspace{40pt}  \bullet\hspace{5pt} & D_{\pm}:=-i\,d/dx,\;
  \mathrm{Dom}(D_{\pm})=H_{\pm}^{1},\;
  D_{\pm}^{2m}:=|D_{\pm}|^{2m},\;
  \mathrm{Dom}(D_{\pm}^{2m})=H_{\pm}^{2m};\hspace{20pt} \\
  \bullet\hspace{5pt}& V(x)=\sum_{k\in \mathbb{Z}}\widehat{V}(2k)e^{i\,2k\pi x}\in
  \mathcal{D}'_{1}; \\
  \bullet\hspace{5pt}& u\in \mathrm{Dom}(S_{\pm}).
\end{align*}
Here by the $H_{\pm}^{1}\equiv H_{\pm}^{1}[0,1]$ and
$H_{\pm}^{2m}\equiv H_{\pm}^{2m}[0,1]$ we denote the Sobolev
spaces of functions that are 1-periodic and 1-semiperiodic on the
interval $[0,1]$, and $\mathcal{D}'_{1}$ denotes the space of
1-periodic distributions \cite[p. 115]{Vld}.

In this paper, we give sufficient conditions for the operators
$S_{\pm}(V)$ to exist as form-sums, conduct a detailed study of
their \textit{qualitative} properties, prove theorems about their
\textit{approximation} and \textit{spectrum decomposition}. The
approximation theorem gives another definition of the operators
$S_{\pm}(V)$ as a limit, in the generalized convergence sense
\cite[Ch. IV, \S 2.6]{Kt}, of a sequence of operators with smooth
potentials.

Earlier in \cite{Mlb, MiMo1, MiMo2}, the authors have carried out
a detailed study of the differential operators $L_{\pm}(V)$
generated on the finite interval by the same differential
expressions as the operators $S_{\pm}(V)$ but defined on the
\textit{negative} Sobolev spaces $H_{\pm}^{-m}$. The case $m=1$
for operators $L_{\pm}(V)$ was treated in \cite{KpMh, Mhr} (see
also closely related papers \cite{MiSb, HrMy, Kor, DjMi}).

So, for an arbitrary $s\in \mathbb{R}$, the Sobolev spaces of
1-periodic and 1-semiperiodic functions or distributions are
defined in a natural fashion by means of their Fourier
coefficients,
\begin{align*}
  H_{+}^{s}\equiv H_{+}^{s}[0,1] & :=\left\{f=\sum_{k\in\mathbb{Z}}\widehat{f}(2k)e^{i 2k\pi
  x}\left|\;\parallel f\parallel_{H_{+}^{s}}<\infty\right.\right\}, \\
  \parallel f\parallel_{H_{+}^{s}}&:=\left(\sum_{k\in\mathbb{Z}}
  \langle 2k\rangle^{2s}\mid\widehat{f}(2k)\mid
  ^{2}\right)^{1/2},\quad \langle k\rangle:=1+|k|, \\
  \widehat{f}(2k)&:=\langle f,e^{i 2k\pi x}\rangle_{+}, \quad
  k\in\mathbb{Z};
\end{align*}
and
\begin{align*}
  H_{-}^{s}\equiv H_{-}^{s}[0,1] & :=\left\{f=\sum_{k\in\mathbb{Z}}\widehat{f}(2k+1)e^{i (2k+1)\pi
  x}\left|\;\parallel f\parallel_{H_{-}^{s}}<\infty\right.\right\}, \\
  \parallel f\parallel_{H_{-}^{s}}&:=\left(\sum_{k\in\mathbb{Z}}
  \langle 2k+1\rangle^{2s}\mid\widehat{f}(2k+1)\mid
  ^{2}\right)^{1/2},\quad \langle k\rangle=1+|k|, \\
  \widehat{f}(2k+1)&:=\langle f,e^{i (2k+1)\pi x}\rangle_{-}, \quad
  k\in\mathbb{Z}.
\end{align*}
By $\langle\cdot,\cdot\rangle_{+}$ and
$\langle\cdot,\cdot\rangle_{-}$ we denote the sesquilinear forms
that define the pairing between the dual spaces $H_{\pm}^{s}$ and
$H_{\pm}^{-s}$ with respect to the zero space $L_{2}(0,1)$; these
pairings are obtained by extending the inner product in
$L_{2}(0,1)$ by continuity \cite[p. 47]{Be},

\begin{equation*}
  (f,g)=\int_{0}^{1}f(x)\overline{g(x)}\,dx,\qquad f,g\in
  L_{2}(0,1).
\end{equation*}

It will be useful to notice that the two-sided scales of Sobolev
spaces $\{H_{\pm}^{s}\}_{s\in \mathbb{R}}$ coincide up to
equivalent norms with scales generated by powers of the
non-negative selfadjoint operators $|D_{\pm}|$ \cite[Ch. II, \S
2.1]{Gor}.

The Sobolev spaces
\begin{equation*}
  H_{per}^{s}\equiv H_{per}^{s}[-1,1]
  :=\left\{f=\sum_{k\in\mathbb{Z}}\widehat{f}(k)e^{i k\pi
  x}\left|\;\parallel f\parallel_{H_{per}^{s}}<\infty\right.\right\},\qquad s\in
  \mathbb{R},
\end{equation*}
of 2-periodic elements (functions or distributions) are defined in
a similar way.

Now, we are ready to formulate the main results obtained in the
paper. But first recall that an operator $A$ on a Hilbert space is
said to be $m$-\textit{sectorial} if its numerical range
$\Theta(A)$, i.e., the set
\begin{equation*}
  \Theta(A):=(A u,u),\quad u\in \mathrm{Dom}(A),\quad \|u\|=1,
\end{equation*}
is contained in a sector of the complex plane,
\begin{align*}
  \Theta(A)&\subseteq \mathrm{Sect}(\gamma,\theta), \\
  \mathrm{Sect}(\gamma,\theta)&:=\{\lambda\in \mathbb{C}\left|\;\right.
  |\arg(\lambda-\gamma)|\leq \theta\},\qquad
  0\leq\theta<\frac{\pi}{2},
\end{align*}
and the exterior of the sector $\mathrm{Sect}(\gamma,\theta)$
belongs to the resolvent set $\mathrm{Resol}(A)$ of the operator
$A$ \cite[Ch. V, \S 3.10]{Kt}.
\begin{theoremEng}\label{ExtTh}
Let a 1-periodic complex-valued distribution V(x) be in the space
$H_{+}^{-m}$. Then the operators $S_{\pm}(V)$ are well defined on
the Hilbert space $L_{2}(0,1)$ as $m$-sectorial operators ---
form-sums,
\begin{equation*}
  S_{\pm}(V)=D_{\pm}^{2m}\dotplus V(x),
\end{equation*}
associated with densely defined, closed, sectorial sesquilinear
forms defined on $L_{2}(0,1)$ by
\begin{equation*}
  t_{S_{\pm}}[u,v]\equiv t_{\pm}[u,v]:=\langle
  D_{\pm}^{2m}u,v\rangle_{\pm}+\langle
  V(x) u,v\rangle_{\pm},\quad \mathrm{Dom}(t_{S_{\pm}})=H_{\pm}^{m},
\end{equation*}
and act on the dense domains
\begin{equation*}
  \mathrm{Dom}(S_{\pm})=\left\{u\in
  H_{\pm}^{m}\,\left|\,D_{\pm}^{2m}u+V(x) u
  \in L_{2}(0,1)\right.\right\}
\end{equation*}
as
\begin{equation*}
  S_{\pm}(V)u=D_{\pm}^{2m}u+V(x)u,\qquad u\in \mathrm{Dom}(S_{\pm}).
\end{equation*}
\end{theoremEng}

Let us remark that, in virtue of the convolution lemma (see Lemma
\ref{CnvLm} below), a 1-periodic complex-valued distribution
$V(x)\in H_{+}^{-m}$ defines, on the Hilbert space $L_{2}(0,1)$,
two sesquilinear forms,
\begin{align*}
  t_{V}^{+}[u,v] & :=\langle V(x)\cdot u,v\rangle_{+},\quad u,v\in H_{+}^{m}, \\
  t_{V}^{-}[u,v] & :=\langle V(x)\cdot u,v\rangle_{-},\quad u,v\in H_{-}^{m},
\end{align*}
where $V(x)\cdot u$ denotes the formal product, which converges in
the Sobolev spaces $H_{\pm}^{-m}$, of the Fourier series of the
distribution $V(x)\in H_{+}^{-m}$ and the function $u\in
H_{\pm}^{m}$.

If the distribution $V(x)$ has additional smoothness in the scale
$\{H_{\pm}^{s}\}_{s\in \mathbb{R}}$ of the Hilbert spaces, then
functions in the domains of the operators $S_{\pm}(V)$ have an
additional regularity.
\begin{theoremEng}\label{RglTh}
Let $V(x)\in H_{+}^{-m\alpha}$, $\alpha\in [0,1]$. Then the
inclusion
\begin{equation*}
\mathrm{Dom}(S_{\pm})\subseteq H_{\pm}^{m(2-\alpha)}
\end{equation*}
holds.
\end{theoremEng}

In the case $\alpha\neq 0$, i.e., for
\begin{equation*}
  V(x)\in H_{+}^{-m\alpha},\quad \alpha\in (0,1],
\end{equation*}
the question about locality of the operators $S_{\pm}(V)$ is
meaningful. Let us recall that an operator $A$ on a function space
is called \textit{local} if
\begin{equation*}
  \mathrm{supp}(A u)\subseteq \mathrm{supp}(u),\qquad u\in
  \mathrm{Dom}(A).
\end{equation*}
For the Hilbert space $L_{2}(0,1)$, this is equivalent to the
following:
\begin{equation*}
  u\left|_{(\alpha,\beta)}\right.=0\quad \Rightarrow
  A u\left|_{(\alpha,\beta)}\right.=0,\qquad u\in
  \mathrm{Dom}(A),\quad(\alpha,\beta)\subset [0,1].
\end{equation*}
\begin{theoremEng}\label{LclTh}
If $V(x)\in H_{+}^{-m}$, the operators $S_{+}(V)$ and $S_{-}(V)$
are local.
\end{theoremEng}

The following theorem describes qualitative properties of the
operators $S_{\pm}(V)$.
\begin{theoremEng}\label{PrpTh}
Let a 1-periodic complex-valued distribution $V(x)$ be in the space
$H_{+}^{-m}$.
\begin{itemize}
  \item [(a)] The operators $S_{\pm}(V)$ are $m$-sectorial with
  respect to an arbitrary angle containing the positive half-axis.
  \item [(b)] The operators $S_{\pm}(V)$ are selfadjoint if and
  only if the distribution $V(x)$ is real-valued, i.e., if
  \begin{equation*}
  \widehat{V}(2k)=\overline{\widehat{V}(-2k)},\quad k\in
  \mathbb{Z}.
  \end{equation*}
  \item [(c)] The operators $S_{\pm}(V)$ have discrete spectra.
\end{itemize}
\end{theoremEng}

The following Theorem \ref{AppTh} allows to give another
alternative definition of the operators $S_{\pm}(V)$ described in
Theorem \ref{ExtTh}.
\begin{theoremEng}\label{AppTh}
Let $V_{n}(x)$, $n\in \mathbb{N}$, and $V(x)$ be defined on the
space $H_{+}^{-m}$, and suppose that
\begin{equation*}
  V_{n}(x)\overset{H_{+}^{-m}}\longrightarrow V(x),\quad
  n\rightarrow\infty.
\end{equation*}
Then the operators $S_{\pm}^{(n)}\equiv S_{\pm}(V_{n})$ converge
to the operators $S_{\pm}\equiv S_{\pm}(V)$ in the uniform
resolvent convergent sense,
\begin{equation*}
  \|R(\lambda,S_{\pm}^{(n)})-R(\lambda,S_{\pm})\|\rightarrow 0,
  \quad
  n\rightarrow\infty.
\end{equation*}
\end{theoremEng}

So, by virtue of Theorem \ref{AppTh}, the operators $S_{\pm}(V)$ can
be defined as a limit of a sequence of the operators $S_{\pm}^{(n)}$
with smooth potentials $V_{n}(x)$ in the generalized convergence
sense \cite[Ch.~IV, \S 2.6]{Kt}.

As an example, consider
\begin{equation*}
  V(x)=\sum_{k\in \mathbb{Z}}\widehat{V}(2k)e^{i\,2k\pi x}\in
  H_{+}^{-m},
\end{equation*}
the trigonometric polynomials
\begin{equation*}
  V_{n}(x)=\sum_{|k|\leq n}\widehat{V}(2k)e^{i\,2k\pi x}\in
  H_{+}^{\infty},
\end{equation*}
form the necessary sequence,
\begin{equation*}
  V_{n}(x)\overset{H_{+}^{-m}}\longrightarrow V(x),\quad
  n\rightarrow\infty,
\end{equation*}
which yields the convergence
\begin{equation*}
  \|R(\lambda,S_{\pm}^{(n)})-R(\lambda,S_{\pm})\|\rightarrow 0,\quad
  n\rightarrow\infty.
\end{equation*}

Due to Theorem \ref{AppTh} we also have that
\begin{equation*}
  \sigma(S_{\pm}^{(n)})\rightarrow \sigma(S_{\pm}),\quad n\rightarrow\infty
\end{equation*}
where the convergence of spectra is upper semicontinuous in general
\cite[Ch. IV, \S 3.1]{Kt} and, for real-valued potentials, it is
continuous \cite[Theorem VIII.23 and Theorem VIII.24]{ReSi1}; by
$\sigma(S_{\pm}^{(n)})$ and $\sigma(S_{\pm})$ we denote unordered
spectra of the corresponding operators.

Now, let us consider, on the Hilbert space $L_{2}(-1,1)$, the
$m$-sectorial operators --- form-sums $S(V)$ with 1-periodic
complex-valued potentials that are distributions $V(x)\in
H_{per}^{-m}$, i.e., $\widehat{V}(2k+1)=0\; \forall k\in
\mathbb{Z}$,
\begin{align*}
L_{2}(-1,1):\hspace{40pt}   & S\equiv S(V):=D^{2m}\dotplus V(x), \\
m\in \mathbb{N}\hspace{45pt}  \bullet\hspace{5pt} & D:=-i\,d/dx,\;
  \mathrm{Dom}(D)=H_{per}^{1},\;
  D^{2m}:=|D|^{2m},\;
  \mathrm{Dom}(D^{2m})=H_{per}^{2m};\hspace{20pt} \\
  \bullet\hspace{5pt}& V(x)=\sum_{k\in \mathbb{Z}}\widehat{V}(2k)e^{i\,2k\pi x}\in
  H_{per}^{-m}; \\
  \bullet\hspace{5pt}& \mathrm{Dom}(S)=\left\{u\in
  H_{per}^{m}\,\left|\,D^{2m}u+V(x) u
  \in L_{2}(-1,1)\right.\right\}.
\end{align*}
Analogs of Theorems \ref{RglTh}, \ref{LclTh}, \ref{PrpTh} and
\ref{AppTh} hold for the operators $S(V)$. In particular, they
have discrete spectra.

Let us study the structure of spectra of the operators $S(V)$,
$S_{+}(V)$, and $S_{-}(V)$ in more details.

Denote by $\mathrm{spec}(A)$ the discrete spectrum of the operator
$A$, taking into account the algebraic multiplicity of the
eigenvalues that ordered lexicographically. Namely, we will say
that an eigenvalue $\lambda_{k}$ precedes an eigenvalue
$\lambda_{k+1}$ for $k\in \mathbb{Z}_{+}$ if
\begin{equation*}
  \mathrm{Re}\,\lambda_{k}<\mathrm{Re}\,\lambda_{k+1},\quad\mbox{or}\quad
  \mathrm{Re}\,\lambda_{k}=\mathrm{Re}\,\lambda_{k+1}\quad\mbox{and}
  \quad\mathrm{Im}\,\lambda_{k}\leq
  \mathrm{Im}\,\lambda_{k+1},\qquad k\in \mathbb{Z}_{+}.
\end{equation*}

It is easy to see that
\begin{align*}
  \mathrm{spec}(D^{2m}) & =\{0; 1,1; 2^{2m},2^{2m};
  \ldots; (2k-1)^{2m},(2k-1)^{2m};
  (2k)^{2m},(2k)^{2m};\ldots\}\cdot \pi^{2m}, \\
  \mathrm{spec}(D_{+}^{2m}) & =\{0;2^{2m},2^{2m};
  \ldots;(2k)^{2m},(2k)^{2m};\ldots\}\cdot \pi^{2m},
  \\
  \mathrm{spec}(D_{-}^{2m}) &
  =\{1,1;3^{2m},3^{2m};\ldots;
  (2k-1)^{2m},(2k-1)^{2m};\ldots\}\cdot \pi^{2m}.
\end{align*}
And thus we get
\begin{equation*}
  \mathrm{spec}(D^{2m})=\mathrm{spec}(D_{+}^{2m})\sqcup
  \mathrm{spec}(D_{-}^{2m})\qquad\text{(the disjoint sum)}.
\end{equation*}
The following Theorem \ref{DcmTh} about spectra decomposition is a
non-trivial generalization of the last equality for the perturbed
m-sectorial operators --- form-sums $S(V)$, $S_{+}(V)$, and
$S_{-}(V)$.
\begin{theoremEng}\label{DcmTh}
Let $S(V)$, $S_{+}(V)$ and $S_{-}(V)$ be the m-sectorial
operators, where the potential $V(x)$ is a 1-periodic
complex-valued distribution from the Sobolev spaces $H_{per}^{-m}$
and $H_{+}^{-m}$ for the first and the second two operators,
respectively. Then
\begin{equation*}
  S(V)=S_{+}(V)\oplus S_{-}(V),
\end{equation*}
and we have the decomposition
\begin{equation*}
  \mathrm{spec}(S)=\mathrm{spec}(S_{+})\cup
  \mathrm{spec}(S_{-}).
\end{equation*}
\end{theoremEng}

A part of results are announced in \cite{MiMo4} and contained in
\cite{MiMo3}.
\section{The proofs}
At first, we will recall some known facts and results that will be
necessary.

Consider the Hilbert spaces of two-sided weighted sequences,
\begin{align*}
  h^{s}&\equiv h^{s}(\mathbb{Z};\mathbb{C}),\quad s\in \mathbb{R},
  \\
  h^{s}&:=\left\{a=\left(a(k)\right)_{k\in \mathbb{Z}}\,\left|\,
  \|a\|_{h^{s}}<\infty\right.\right\}, \\
  (a,b)_{h^{s}}&:=\sum_{k\in\mathbb{Z}}
  \langle k\rangle^{2s}a(k)\,\overline{b(k)},\quad \langle
  k\rangle=1+|k|, \\
  \|a\|_{h^{s}}&:=\left(\sum_{k\in\mathbb{Z}}
  \langle k\rangle^{2s}|a(k)| ^{2}\right)^{1/2}.
\end{align*}
The Fourier transform establishes an isometric isomorphisms between
the Sobolev spaces $H_{per}^{s}$, $H_{\pm}^{s}$ and the Hilbert
spaces $h^{s}$ of two-sided weighted sequences,
\begin{align*}
  \mathcal{F}: H_{per}^{s}\ni f & \mapsto (\widehat{f})=
  \left(\widehat{f}(k)\right)_{k\in \mathbb{Z}}\in h^{s}, \\
  \mathcal{F}_{+}: H_{+}^{s}\ni f & \mapsto (\widehat{f})=
  \left(\widehat{f}(2k)\right)_{k\in \mathbb{Z}}\in h^{s}, \\
  \mathcal{F}_{-}: H_{-}^{s}\ni f & \mapsto (\widehat{f})=
  \left(\widehat{f}(2k+1)\right)_{k\in \mathbb{Z}}\in h^{s}.
\end{align*}
This, together with the convolution lemma (see bellow), allows to
give sufficient conditions of existence of the formal product
\begin{equation*}
  V(x)\cdot u(x)=\sum_{k\in \mathbb{Z}}\sum_{j\in \mathbb{Z}}
  \widehat{V}(k-j)\widehat{u}(j)\,e^{i\,k\pi x}.
\end{equation*}
To this end, introduce in the scale of the Hilbert spaces of
two-sided weighted sequences $\{h^{s}\}_{s\in \mathbb{R}}$ a
commutative convolution operation. For arbitrary sequences
\begin{equation*}
  a=\left(a(k)\right)_{k\in\mathbb{Z}}\quad\text{and}\quad
  b=\left(b(k)\right)_{k\in\mathbb{Z}},
\end{equation*}
it is defined in a natural fashion,
\begin{align*}
  (a,b)&\mapsto a*b, \\
  (a*b)(k)&:=\sum_{j\in\mathbb{Z}}a(k-j)b(j).
\end{align*}

The following known lemma (see, for example \cite[Lemma
1.5.4]{Mhr}) holds.
\begin{lemmaEng}[Convolution lemma]\label{CnvLm}
Let $s,r\geq 0$, and $t\in \mathbb{R}$ with $t\leq \min(s,r)$.
\begin{itemize}
  \item [(\verb"I")] If $s+r-t>1/2$, then the convolution map
  \begin{equation*}
    (a,b)\mapsto a*b
  \end{equation*}
  is continuous when viewed as the maps
  \begin{equation*}
    (a)\; h^{r}\times h^{s}\rightarrow h^{t},\hspace{100pt}
    (b)\;h^{-t}\times h^{s}\rightarrow h^{-r}.
  \end{equation*}
  \item [(\verb"II")] If $s+r-t<1/2$, then this statement fails to
  hold.
\end{itemize}
\end{lemmaEng}

\subsection{Proof of Theorem \ref{ExtTh}}
Due to the convolution lemma for
\begin{equation*}
  V(x)\in H_{+}^{-m}\quad\text{and}\quad u(x)\in H_{\pm}^{m},
\end{equation*}
the products $V(x)\cdot u(x)$ are well defined in the Sobolev
spaces $H_{\pm}^{-m}$. Therefore, the sesquilinear forms
\begin{align*}
  t_{V}^{+}[u,v]&=\langle V(x)u,v \rangle_{+},\qquad
  \mathrm{Dom}(t_{V}^{+})=H_{+}^{m}, \\
  t_{V}^{-}[u,v]&=\langle V(x)u,v \rangle_{-},\qquad
  \mathrm{Dom}(t_{V}^{-})=H_{-}^{m}
\end{align*}
are well defined in the Hilbert space $L_{2}(0,1)$.

Further, set
\begin{align*}
  \tau_{+}[u,v]&:=\langle D_{+}^{2m}u,v \rangle_{+},\qquad
  \mathrm{Dom}(\tau_{+})=H_{+}^{m}, \\
  \tau_{-}[u,v]&:=\langle D_{-}^{2m}u,v \rangle_{-},\qquad
  \mathrm{Dom}(\tau_{-})=H_{-}^{m}.
\end{align*}
The sesquilinear forms $\tau_{\pm}[u,v]$ are well defined in the
Hilbert space $L_{2}(0,1)$, they are densely defined, closed, and
nonnegative.

The following fundamental assertion is true.
\begin{propositionEng}\label{RltBndPrp}
The sesquilinear forms $t_{V}^{\pm}[u,v]$ are $\tau_{\pm}$-bounded
with $\tau_{\pm}$-boundary that equals zero, i.e., we have
$V(x)\prec\prec D_{\pm}^{2m}$.
\end{propositionEng}
\begin{proof}
Represent the 1-periodic distribution
\begin{equation*}
  V(x)=\sum_{k\in \mathbb{Z}}\widehat{V}(2k)e^{i\,2k\pi x}\in
  H_{+}^{-m}
\end{equation*}
as the sum
\begin{equation}\label{eq_10}
  V(x)=V_{0}(x)+V_{\delta}(x),
\end{equation}
where $V_{0}(x)$ is a smooth function and $V_{\delta}(x)$ is a
distribution with an arbitrarily small norm,
\begin{equation*}
  V_{0}(x)\in H_{+}^{m},
\end{equation*}
and
\begin{equation*}
  V_{\delta}(x)\in H_{+}^{-m}\quad\text{with}\quad
  \|V_{\delta}\|_{H_{+}^{-m}}\leq \frac{\delta}{C_{m}}.
\end{equation*}
The constant $C_{m}$ is defined from the convolution lemma and is
fixed. The decomposition \eqref{eq_10} is possible, since
$H_{+}^{m}$ is densely embedded into the space $H_{+}^{-m}$.

So, for
\begin{equation*}
  u\in\mathrm{Dom}(\tau_{\pm})\subset\mathrm{Dom}(t_{V}^{\pm}),
\end{equation*}
we have
\begin{align*}
  \left|t_{V}^{\pm}[u]\right|=\left|\langle V(x) u,u\rangle_{\pm}\right|
  &\leq \left|\langle V_{0}(x) u,u\rangle_{\pm}\right|
  +\left|\langle V_{\delta}(x) u,u\rangle_{\pm}\right| \\
  & \leq \left\|V_{0}(x) u\right\|_{L_{2}(0,1)}\left\|u\right\|_{L_{2}(0,1)}
  +\left\|V_{\delta}(x) u\right\|_{H_{\pm}^{-m}}\left\|u\right\|_{H_{\pm}^{m}} \\
  & \leq C_{m}\left\|V_{0}(x)\right\|_{H_{+}^{m}}\left\|u\right\|_{L_{2}(0,1)}^{2}
  +\delta\left\|u\right\|_{H_{\pm}^{m}}^{2}.
\end{align*}
Taking to account that
\begin{equation*}
  \left\|u\right\|_{H_{\pm}^{m}}^{2}\leq \left\|u\right\|_{L_{2}(0,1)}^{2}
  +\|u^{(m)}\|_{L_{2}(0,1)}^{2}
  =\left\|u\right\|_{L_{2}(0,1)}^{2}+\langle D_{\pm}^{2m}u,u\rangle_{\pm}
\end{equation*}
for an arbitrary $\delta>0$ we obtain the necessary estimate,
\begin{equation}\label{eq_12}
  \left|t_{V}^{\pm}[u]\right|\leq \delta \tau_{\pm}[u]
  +\left(C_{m}\left\|V_{0}\right\|_{H_{+}^{m}[0,1]}+\delta\right)\|u\|_{L_{2}(0,1)}^{2}.
\end{equation}

The proof is complete.
\end{proof}
Proposition \ref{RltBndPrp}, together with \cite[Theorem
IV.1.33]{Kt} yields the following
\begin{corollaryEng}\label{crl_10}
The sesquilinear forms
\begin{equation*}
  t_{\pm}[u,v]:=\langle D_{\pm}^{2m}u,v\rangle_{\pm}+\langle V(x)u,v\rangle_{\pm},
  \quad \mathrm{Dom}(t_{\pm})=H_{\pm}^{m}
\end{equation*}
are densely defined, closed, and sectorial in the Hilbert space
$L_{2}(0,1)$.
\end{corollaryEng}

According to the first representation theorem \cite[Theorem
VI.2.1]{Kt}, there exist $m$-sectorial operators $S_{\pm}(V)$
associated with the forms $t_{\pm}[u,v]$ such that
\begin{itemize}
  \item [i)] $\mathrm{Dom}(S_{\pm})\subseteq
  \mathrm{Dom}(t_{\pm})$ and
  \begin{equation*}
    t_{\pm}[u,v]=(S_{\pm}u,v)
  \end{equation*}
  for every  $u\in \mathrm{Dom}(S_{\pm})$ and $v\in
  \mathrm{Dom}(t_{\pm})$;
  \item [ii)] $\mathrm{Dom}(S_{\pm})$ are cores of $t_{\pm}[u,v]$;
  \item [iii)] if $u\in \mathrm{Dom}(t_{\pm})$, $w\in L_{2}(0,1)$,
  and
  \begin{equation*}
    t_{\pm}[u,v]=(w,v)
  \end{equation*}
  holds for every $v$ belonging to the cores of $t_{\pm}[u,v]$, then $u\in
  \mathrm{Dom}(S_{\pm})$ and $S_{\pm}(V)u=w$. The $m$-sectorial
  operators $S_{\pm}(V)$ are uniquely defined by condition
  i).
\end{itemize}

Now, investigate the operators $S_{\pm}(V)$ associated with the
forms $t_{\pm}[u,v]$ in more details.

Let
\begin{equation*}
  u\in \mathrm{Dom}(S_{\pm})\quad\text{and}\quad v\in
  \mathrm{Dom}(t_{\pm}).
\end{equation*}
Then we have
\begin{align*}
 t_{\pm}[u,v]&=\langle D_{\pm}^{2m}u,v\rangle_{\pm}+\langle V(x)u,v\rangle_{\pm}
 =\langle D_{\pm}^{2m}u+V(x)u,v\rangle_{\pm}
=\left(S_{\pm}u,v\right)=\langle S_{\pm}u,v\rangle_{\pm}.
\end{align*}
This shows that we have the equality
\begin{equation*}
  \langle D_{\pm}^{2m}u+V(x) u,v\rangle_{\pm}=\langle S_{\pm}u,v\rangle_{\pm},\quad v\in
  H_{\pm}^{m},
\end{equation*}
of linear forms. So, we can conclude that
\begin{equation*}
    S_{\pm}(V)u=D_{\pm}^{2m}u+V(x) u\in L_{2}(0,1),\quad u\in
    \mathrm{Dom}(S_{\pm}),
\end{equation*}
and that the inclusions
\begin{equation*}
  \mathrm{Dom}(S_{\pm})\subseteq
  \left\{u\in H_{\pm}^{m}\,\left|\,D_{\pm}^{2m}u+V(x) u\in L_{2}(0,1)\right.\right\}
\end{equation*}
hold. It remains to verify that the inverse inclusions hold.

Let
\begin{equation*}
  u\in \left\{u\in H_{\pm}^{m}\,\left|\,D_{\pm}^{2m}u+V(x) u\in L_{2}(0,1)\right.\right\}
  \quad\text{and}\quad v\in \mathrm{Dom}(t_{\pm}).
\end{equation*}
Then
\begin{align*}
 t_{\pm}[u,v]&=\langle D_{\pm}^{2m}u,v\rangle_{\pm}+
 \langle V(x) u,v\rangle_{\pm}=\langle D_{\pm}^{2m}u+V(x)
 u,v\rangle_{\pm}=\left(D_{\pm}^{2m}u+V(x) u,v\right),
\end{align*}
and using the first representation theorem iii) (see above) we get
the necessary estimate,
\begin{equation*}
  u\in \mathrm{Dom}(S_{\pm}),
\end{equation*}
which implies that
\begin{equation*}
  \left\{u\in H_{\pm}^{m}\,\left|\,D_{\pm}^{2m}u+V(x) u
  \in L_{2}(0,1)\right.\right\}\subseteq\mathrm{Dom}(S_{\pm})
\end{equation*}
and
\begin{equation*}
  S_{\pm}(V)u=D_{\pm}^{2m}u+V(x) u\in L_{2}(0,1).
\end{equation*}

So,
\begin{equation*}
  \mathrm{Dom}(S_{\pm})=\left\{u\in H_{\pm}^{m}\,\left|\,D_{\pm}^{2m}u+V(x) u
  \in L_{2}(0,1)\right.\right\}
\end{equation*}
and
\begin{equation*}
  S_{\pm}(V)u=D_{\pm}^{2m}u+V(x) u\in L_{2}(0,1),\qquad u\in
  \mathrm{Dom}(S_{\pm}).
\end{equation*}

Theorem \ref{ExtTh} is proved completely.
\begin{remarkEng}\label{rmk_10}
Throughout the rest of the paper we will often use the notations
\begin{equation*}
  t_{S_{\pm}}[u,v]\equiv t_{\pm}[u,v]
\end{equation*}
to underline the dual relations between the sesquilinear forms
$t_{\pm}[u,v]$ and the associated with them operators
$S_{\pm}(V)$, see \cite[Theorem VI.2.7]{Kt}.
\end{remarkEng}

\subsection{Proof of Theorem \ref{RglTh}}
Let the 1-periodic distribution $V(x)$ belong to the space
$H_{+}^{-m\alpha}$, $\alpha\in [0,1]$. Then for any $u\in
\mathrm{Dom}(S_{\pm})$, due to the convolution lemma, we have
\begin{equation*}
  V(x)u\in H_{+}^{-m\alpha},
\end{equation*}
and therefore
\begin{equation*}
  D_{\pm}^{2m}u\in H_{\pm}^{-m\alpha}.
\end{equation*}
From this we conclude that
\begin{equation*}
  u\in H_{\pm}^{m(2-\alpha)}.
\end{equation*}

\subsection{Proof of Theorem \ref{LclTh}}
Let
\begin{equation*}
  u\in \mathrm{Dom}(S_{\pm})
\end{equation*}
and
\begin{equation*}
  u\left|_{(\alpha,\beta)}\right.=0\quad\text{with} \quad(\alpha,\beta)\subset
  [0,1],
\end{equation*}
and let
\begin{equation*}
  \varphi(x)\in C_{0}^{\infty}[0,1]\quad\text{with} \quad
  \mathrm{supp}(\varphi)\Subset (\alpha,\beta).
\end{equation*}
Then we have
\begin{align*}
 (S_{\pm}u,\varphi)&=\langle S_{\pm}u,\varphi\rangle_{\pm}
 =\langle D_{\pm}^{2m}u+V(x) u,\varphi\rangle_{\pm}
 =\langle D_{\pm}^{2m}u,\varphi\rangle_{\pm}+\langle V(x)
 u,\varphi\rangle_{\pm} \\
 &=\langle u,D_{\pm}^{2m}\varphi\rangle_{\pm}+
 \langle V(x),\overline{u}\varphi\rangle_{\pm}
 =\langle V(x),0\rangle_{\pm}=0,
\end{align*}
which yields the necessary statement,
\begin{equation*}
  \left(S_{\pm}u\right)\left|_{(\alpha,\beta)}\right.=0.
\end{equation*}

\subsection{Proof of Theorem \ref{PrpTh}}
(a) The $m$-sectoriality of the operators $S_{\pm}(V)$ have been
proved in Theorem \ref{ExtTh}. Let us prove the second part of the
assertion, i.e., we need to show that for any $\varepsilon>0$ and
some constant $c_{\varepsilon}\geq 0$ the following estimates
hold:
\begin{equation*}
  \left|\arg\left((S_{\pm}+c_{\varepsilon}Id)u,u\right)\right|\leq
  \varepsilon,\quad u\in \mathrm{Dom}(S_{\pm}).
\end{equation*}
For this we have to make sure that
\begin{equation*}
  \left|\mathrm{Im}(S_{\pm}u,u)\right|
  \leq \varepsilon\mathrm{Re}(S_{\pm}u,u)+c_{\varepsilon}\|u\|_{L_{2}(0,1)}^{2},
  \quad u\in\mathrm{Dom}(S_{\pm}),
\end{equation*}
for any $\varepsilon>0$ and some constant $c_{\varepsilon}\geq 0$.

So, take $0<\varepsilon<1/2$. From Proposition \ref{RltBndPrp}
(see \eqref{eq_12}) we get
\begin{equation*}
  \left|t_{V}^{\pm}[u]\right|\leq  \frac{\varepsilon}{2}\tau_{\pm}[u]
  +\left(C_{m}\left\|V_{0}(x)\right\|_{H_{+}^{m}}
  +\frac{\varepsilon}{2}\right)\|u\|_{L_{2}(0,1)}^{2},
  \quad u\in \mathrm{Dom}(\tau_{\pm}),
\end{equation*}
and, hence,
\begin{equation*}
  -\varepsilon\mathrm{Re}\,t_{V}^{\pm}[u]\leq  \frac{\varepsilon}{2}\tau_{\pm}[u]
  +\left(C_{m}\left\|V_{0}(x)\right\|_{H_{+}^{m}}
  +\frac{\varepsilon}{2}\right)\|u\|_{L_{2}(0,1)}^{2},
  \quad u\in \mathrm{Dom}(\tau_{\pm}).
\end{equation*}
Further, taking into account that
\begin{align*}
  \mathrm{Re}(S_{\pm}u,u) & =\langle D_{\pm}^{2m}u,u\rangle_{\pm}
  +\mathrm{Re}\langle V(x) u,u\rangle_{\pm}, \\
  \mathrm{Im}(S_{\pm}u,u) & =\mathrm{Im}\langle V(x) u,u\rangle_{\pm},
  \quad u\in \mathrm{Dom}(S_{\pm})
\end{align*}
we obtain the necessary estimates,
\begin{align*}
  \left|\mathrm{Im}(S_{\pm}u,u)\right| & \leq \left|\langle V(x) u,u\rangle_{\pm}\right|\leq
  \frac{\varepsilon}{2}\tau_{\pm}[u]+\left(C_{m}\left\|V_{0}(x)\right\|_{H_{+}^{m}}
  +\frac{\varepsilon}{2}\right)\|u\|_{L_{2}(0,1)}^{2} \\
  & \leq \varepsilon\left(\tau_{\pm}[u]+\mathrm{Re}\,t_{V}^{\pm}[u]\right)
  +\left(2C_{m}\left\|V_{0}(x)\right\|_{H_{+}^{m}}
  +\varepsilon\right)\|u\|_{L_{2}(0,1)}^{2} \\
  & =\varepsilon\mathrm{Re}(S_{\pm}u,u)+c_{\varepsilon}\|u\|_{L_{2}(0,1)}^{2},
  \qquad u\in \mathrm{Dom}(S_{\pm}).
\end{align*}

(b) Let the 1-periodic distribution $V(x)$ be real-valued. Then the
sesquilinear forms $t_{S_{\pm}}[u,v]$ are symmetric and,
consequently, in virtue of \cite[Theorem VI.2.7]{Kt} (also see the
KLMN theorem \cite[Theorem X.17]{ReSi2}), the operators are
selfadjoint.

Conversely, let the operators $S_{\pm}(V)$ be selfadjoint. In the
case of a non-real-valued distribution $V(x)$, the operators
$S_{\pm}(V)$ are not symmetric either. This contradiction allows
to make conclusion that the distribution $V(x)$ is real-valued.

(c) From \cite[Theorem VI.3.4]{Kt} and Proposition \ref{RltBndPrp}
we immediately obtain the following.
\begin{propositionEng}\label{CpRsPr}
The resolvent sets of the operators $S_{\pm}(V)$ are non-empty.
Moreover, the resolvents $R(\lambda,S_{\pm}(V))$ of the operators
$S_{\pm}(V)$  are compact.
\end{propositionEng}

Proposition \ref{CpRsPr} implies that the operators $S_{\pm}(V)$
have discrete spectra.

\subsection{Proof of Theorem \ref{AppTh}}
The proof is based on the following proposition.
\begin{propositionEng}\label{GnCnPr}
Let the 1-periodic distributions $V_{n}(x)$, $n\in \mathbb{N}$,
and $V(x)$ be in the Sobolev space $H_{+}^{-m}$. For
\begin{equation*}
  V_{n}(x)\overset{H_{+}^{-m}}\longrightarrow V(x),\quad
  n\rightarrow\infty,
\end{equation*}
the operators
\begin{align*}
  S_{\pm}^{(n)}&\equiv S_{\pm}(V_{n}):=D_{\pm}^{2m}\dotplus V_{n}(x),
  \\
  \mathrm{Dom}(S_{\pm}^{(n)})&=\left\{u\in H_{\pm}^{m}\,\left|\,D_{\pm}^{2m}u+V_{n}(x) u
  \in L_{2}(0,1)\right.\right\},
\end{align*}
converge to the operators
\begin{align*}
  S_{\pm}&\equiv S_{\pm}(V)=D_{\pm}^{2m}\dotplus V(x),
  \\
  \mathrm{Dom}(S_{\pm})&=\left\{u\in H_{\pm}^{m}\,\left|\,D_{\pm}^{2m}u+V(x) u
  \in L_{2}(0,1)\right.\right\},
\end{align*}
in the generalized convergence sense \cite[Ch. IV, \S 2.6]{Kt}.
\end{propositionEng}
\begin{proof}
Set
\begin{equation*}
   t_{S_{\pm}^{(n)}}[u,v]\equiv t_{\pm}^{(n)}[u,v]:=(S_{\pm}^{(n)}u,v),\qquad
   u\in\mathrm{Dom}(S_{\pm}^{(n)}),v\in\mathrm{Dom}(t_{\pm}^{(n)})=H_{\pm}^{m},
\end{equation*}
and recall that
\begin{equation*}
   t_{S_{\pm}}[u,v]\equiv t_{\pm}[u,v]=(S_{\pm}u,v),\qquad
   u\in\mathrm{Dom}(S_{\pm}),v\in\mathrm{Dom}(t_{\pm})=H_{\pm}^{m}.
\end{equation*}
Then, for every $u\in \mathrm{Dom}(t_{\pm})=
\mathrm{Dom}(t_{\pm}^{(n)})=H_{\pm}^{m}$,
\begin{align*}
  \left|t_{\pm}^{(n)}[u]-t_{\pm}[u]\right| &
  =\left|\left\langle (V_{n}(x)-V(x))u,u\right\rangle_{\pm}\right|
  \leq \|(V_{n}(x)-V(x))u\|_{H_{\pm}^{-m}} \|u\|_{H_{\pm}^{m}} \\
  & \leq C_{m}\|V_{n}(x)-V(x)\|_{H_{+}^{-m}}
  \left(\|u\|_{H_{\pm}^{m}}^{2}+\tau_{\pm}[u]\right),
\end{align*}
where the constant $C_{m}$ is defined due to the convolution
lemma, and $\tau_{\pm}[u,v]$, as above,
\begin{equation*}
  \tau_{\pm}[u,v]=\langle D_{\pm}^{2m}u,v\rangle_{\pm},
  \qquad \mathrm{Dom}(\tau_{\pm})=H_{\pm}^{m},
\end{equation*}
are sesquilinear, densely defined, closed and nonnegative forms.
Since the forms
\begin{equation*}
  t_{V}^{\pm}[u,v]=\langle V(x)u,v\rangle_{\pm},
  \qquad \mathrm{Dom}(t_{V}^{\pm})=H_{\pm}^{m},
\end{equation*}
are $\tau_{\pm}$-bonded with zero $\tau_{\pm}$-boundary for an
arbitrary $0<\varepsilon \leq1/2$, the following estimates hold:
\begin{equation*}
  2\left|\mathrm{Re}\,t_{V}^{\pm}[u]\right|
  \leq \tau_{\pm}[u]+ 2\left(C_{m}\left\|V_{0}(x)\right\|_{H_{+}^{m}}+\varepsilon\right)
  \|u\|_{L_{2}(0,1)}^{2},
\end{equation*}
and thus
\begin{equation*}
  2\mathrm{Re}\,t_{V}^{\pm}[u]+\tau_{\pm}[u]+
  2\left(C_{m}\left\|V_{0}(x)\right\|_{H_{+}^{m}}+\varepsilon\right)\|u\|_{L_{2}(0,1)}^{2}
  \geq 0.
\end{equation*}
Taking to account that
\begin{equation*}
  \mathrm{Re}\,t_{\pm}[u]=\tau_{\pm}[u]+\mathrm{Re}\,t_{V}^{\pm}[u]
\end{equation*}
we get the needed estimates,
\begin{align*}
  \left|t_{\pm}^{(n)}[u]-t_{\pm}[u]\right|
  &\leq
  C_{m}\|V_{n}(x)-V(x)\|_{H_{+}^{-m}}\left(\|u\|_{H_{\pm}^{m}}^{2}+\tau_{\pm}[u]\right)
  \leq C_{m}\|V_{n}(x)-V(x)\|_{H_{+}^{-m}} \\
  &\cdot\left(2\mathrm{Re}\,t_{V}^{\pm}[u]+2\tau_{\pm}[u]
  +2\left(C_{m}\left\|V_{0}(x)\right\|_{H_{+}^{m}}+\varepsilon+1/2\right)
  \|u\|_{L_{2}(0,1)}^{2}\right) \\
  &=a_{n}\|u\|_{L_{2}(0,1)}^{2}+b_{n}\mathrm{Re}\,t_{\pm}[u],
\end{align*}
where
\begin{equation*}
  a_{n}=2\left(C_{m}\left\|V_{0}(x)\right\|_{H_{+}^{m}}+1\right)
  \|V_{n}(x)-V(x)\|_{H_{+}^{-m}}\geq 0
\end{equation*}
and
\begin{equation*}
  b_{n}=2C_{m}\|V_{n}(x)-V(x)\|_{H_{+}^{-m}}\geq 0
\end{equation*}
tend to zero as $n\rightarrow\infty$.

To complete the proof it suffices to apply \cite[Theorem
VI.3.6]{Kt}.
\end{proof}

Proposition \ref{GnCnPr} and \cite[Theorem IV.2.25]{Kt}, together
with Proposition \ref{CpRsPr}, give Theorem \ref{AppTh}.

\subsection{Proof of Theorem \ref{DcmTh}}
Let the operators --- form-sums $S(V)$, $S_{+}(V)$, and $S_{-}(V)$
be given with $V(x)$ a 1-periodic complex-valued distribution from
the Sobolev spaces $H_{per}^{-m}$ and $H_{\pm}^{-m}$,
correspondingly.

For an arbitrary $s\in \mathbb{R}$ let us consider the Sobolev
spaces
\begin{align*}
 H_{per}^{s}
 &=\left\{f=\sum_{k\in\mathbb{Z}}\widehat{f}(k)e^{i k\pi
 x}\left|\;\parallel
 f\parallel_{H_{per}^{s}}<\infty\right.\right\}, \\
 H_{+}^{s}
 &=\left\{f=\sum_{k\in\mathbb{Z}}\widehat{f}(2k)e^{i 2k\pi
 x}\left|\;\parallel f\parallel_{H_{+}^{s}}<\infty\right.\right\},
 \\
 H_{-}^{s}
 &=\left\{f=\sum_{k\in\mathbb{Z}}\widehat{f}(2k+1)e^{i (2k+1)\pi
 x}\left|\;\parallel f\parallel_{H_{-}^{s}}<\infty\right.\right\}.
\end{align*}
It should be remarked that
\begin{equation*}
  H_{per}^{0}\equiv L_{2}(-1,1)\quad\text{and}
  \quad H_{+}^{0}\equiv H_{-}^{0}\equiv L_{2}(0,1).
\end{equation*}
Set
\begin{align*}
  H_{per,+}^{s}&:=\left\{f\in H_{per}^{s}\,\left|\,
  \widehat{f}(2k+1)=0\quad \forall k\in \mathbb{Z}\right.\right\},
  \\
  H_{per,-}^{s}&:=\left\{f\in H_{per}^{s}\,\left|\,
  \widehat{f}(2k)=0\quad \forall k\in \mathbb{Z}\right.\right\},
\end{align*}
and thus
\begin{equation*}
  H_{per}^{s}=H_{per,+}^{s}\oplus H_{per,-}^{s},\quad s\in
  \mathbb{R}.
\end{equation*}

Let
\begin{equation*}
  I_{\pm}:\;
  H_{\pm}^{s}\ni f(x)\mapsto f(x)\in H_{per,\pm}^{s},
  \quad s\in \mathbb{R},
\end{equation*}
be extension operators that extend the elements $f(x)\in
H_{\pm}^{s}$ defined on the interval $[0,1]$ to the elements
$f(x)\in H_{per,\pm}^{s}$ defined on the interval $[-1,1]$. The
operators $I_{\pm}$ establish isometric isomorphisms between the
spaces $H_{\pm}^{s}$ and $H_{per,\pm}^{s}$ for $s\in \mathbb{R}$.

Further, let us consider the operators $S(V)$. Since the potentials
$V(x)$ are 1-periodic distributions from the space $H_{per}^{-m}$,
i.e., $V(x)\in H_{per,+}^{-m}$, the operators $S(V)$ are reduced by
the space $H_{per,+}^{-m}$ \cite[Ch. IV, \S 40]{AkhGl}. So, we have
\begin{equation}\label{eq_14}
  S(V)=S_{per,+}(V)\oplus S_{per,-}(V),
\end{equation}
where the operators $S_{per,\pm}(V)$ are defined on the Hilbert
spaces $H_{per,\pm}^{0}$. Taking into account that
\begin{equation*}
  H_{+}^{s}\overset{I_{+}}\simeq H_{per,+}^{s}\quad\text{and}
  \quad H_{-}^{s}\overset{I_{-}}\simeq H_{per,-}^{s}
\end{equation*}
for an arbitrary $s\in \mathbb{R}$ we conclude that the operators
$S_{per,\pm}(V)$ and $S_{\pm}(V)$ are unitary equivalent,
\begin{equation*}
  S_{+}(V)\overset{I_{+}}\simeq S_{per,+}(V)\quad\text{and}
  \quad S_{-}(V)\overset{I_{-}}\simeq S_{per,-}(V).
\end{equation*}
From the latter relations and decomposition \eqref{eq_14}, we
obtain the need statement,
\begin{equation*}
  S(V)=S_{+}(V)\oplus S_{-}(V),
\end{equation*}
which implies
\begin{equation*}
  \mathrm{spec}(S)=\mathrm{spec}(S_{+})\cup
  \mathrm{spec}(S_{-}).
\end{equation*}

The proof of Theorem \ref{DcmTh} is completed.


\end{document}